\def\versiondate{25 Sep.\ 2002}
\input math.macros
\input Ref.macros

\checkdefinedreferencetrue
\continuousfigurenumberingtrue
\theoremcountingtrue
\sectionnumberstrue
\forwardreferencetrue
\citationgenerationtrue
\nobracketcittrue
\hyperstrue
\initialeqmacro

\bibsty{myapalike}
\input\jobname.lbl 

\def\T{{\Bbb T}}
\def\GM{{\ss GM}}  
\def\ip#1{(\changecomma #1)}
\def\bigip#1{\bigl(\bigchangecomma #1\bigr)}

\def\changecomma#1,{#1,\,}
\def\bigchangecomma#1,{#1,\;}
\def\leftchangecomma#1,{#1,\ }
\def\constant#1{{\bf #1}}
\def\trig{{\ss Poly}}  
\def\onp{{\bf p}}  
\def\tilonp{\widetilde{{\bf p}}}  
\def\qonp{{\bf q}}  
\def\bfz{{\bf 0}}
\def\bfo{{\bf 1}}
\def\hl{{\ss HL}}  
\def\supp{{\rm supp}\,}
\def\F{{\cal F}}
\def\G{{\cal G}}

\ifproofmode \relax \else\head{To appear in {\it Geom. Funct. Anal.}}
{Version of \versiondate}\fi 
\vglue20pt

\title{Szeg\H{o} Limit Theorems}

\author{Russell Lyons}

\abstract{%
The first Szeg\H{o} limit theorem has been extended by
Bump-Diaconis and Tracy-Widom to limits of other minors of Toeplitz
matrices.
We use a more geometric method to
extend their results still further. Namely, we allow more general measures and
more general determinants.
We also give a new extension to higher dimensions, which extends a theorem
of \refbauthor{HelLow}.
}

\bottomII{Primary 
47B35. 
Secondary
11C20. 
}
{Toeplitz determinants, minors, spectral factors.}
{Research partially supported by 
NSF grant DMS-0103897.}

\bsection{Introduction}{s.intro}

Let $\lambda$ denote Lebesgue measure on the unit-length circle $\T := \R/\Z$.
For a finite positive measure $\mu$ on $\T$, its Fourier coefficients are
$\widehat\mu(n) := \int_{\T} e^{-2\pi i n t}\,d\mu(t)$.
For $n \ge 0$, define 
$$
D_n(\mu) := \det[\widehat\mu(k-j)]_{0 \le j, k \le n}
\,.
$$
The Szeg\H{o} limit theorems determine the asymptotics of $D_n(\mu)$.
The first Szeg\H{o} limit theorem (\ref b.GZ/, p.~44), which is the one of
concern here, determines $\lim_{n \to\infty} D_{n+1}(\mu)/D_n(\mu)$.
To state this beautiful result of Szeg\H{o}'s in its extended form due to
Kolmogorov and Kre\u{\i}n, define,
for $f \ge 0$ and $\log^+ f \in L^1(\lambda)$, 
the {\bf geometric mean} of $f$ by 
$$
\GM(f) := \exp \int_{\T} \log f \,d\lambda
\,.
$$
According to the arithmetic mean-geometric mean inequality, if $0 \le f \in
L^1(\lambda)$, then $\log^+ f \in L^1(\lambda)$ and $0 \le \GM(f) \le \int f
\,d\lambda$.

\procl t.first \procname{Szeg\H{o} Limit Theorem}
Let $\mu$ be a finite positive measure on $\T$ with infinite support.
Let $f:= [d\mu/d\lambda]$ be the Radon-Nikod\'ym derivative of the
absolutely continuous part of $\mu$.
Then 
$$
\lim_{n \to\infty} D_{n+1}(\mu)/D_n(\mu) = \GM(f)
\,.
\label e.basic
$$
\endprocl

This result has been extended in various ways, three of which we consider
here.
The first two extensions were proved by \ref b.BD/ and \ref b.TW/, while
the third was proved by \ref b.HelLow/.
Our theorems extend theirs further yet.
These extensions concern quotients of determinants where, instead of the
determinant in the numerator having one more particular row and column than
the determinant in the denominator, as in \ref e.basic/, the numerator
contains finitely many more rows and columns arising from inner
products among arbitrary vectors.
In addition, we shall consider the case where $\mu$ is complex, as did
\refbauthor{BD} and \refbauthor{TW}.
Theorems \briefref t.nonneg/ and \briefref t.nonneg-d/ of this paper are
used in \ref b.LyonsSteif:dyn/ for studying entropy and phase transitions
in stationary determinantal processes. 

The approach of \refbauthor{BD} relies on certain identities for symmetric
functions and representations of the symmetric group, while \refbauthor{TW}
proceed via factorization and analysis of Toeplitz operators.
Our approach is different, although it bears some similarities to that of
\refbauthor{TW}. Namely, we go to the level of vectors in Hilbert space
and analyze projections on various subspaces. This more geometric
method appears to be more flexible
and leads to a more transparent proof and formulation of the results.
Although we use row operations that leave determinants unchanged, one could
instead use multivectors and continue a geometric approach. However, to
maintain greater accessibility, we have omitted exterior algebra.

In a brief \ref s.notn/, we recall some general notation and facts from
complex analysis.
This allows us to state and prove our first result in \ref s.positive/, which
concerns the case $\mu \ge 0$.
We turn to the case of complex $\mu$ in \ref s.cplx/.
This requires a lemma about convergence of non-orthogonal linear
projections in Hilbert space, whose proof is relegated to the appendix of
the paper.
\ref b.HelLow/ proved an analogue of \ref t.first/ for higher dimensions.
In \ref s.high/, we extend their theorem by proving analogues of
our results for higher
dimensions; both the real and complex case are treated there.
Actually, \refbauthor{HelLow} stated their result not as a limit of
quotients of determinants, but as an extremum problem, as in \ref b.GZ/.

\bsection{General Notation and Hardy Spaces}{s.notn}

Write $e_n(t) := e^{2\pi i n t}$ and for $f \in L^1(\lambda)$, write $\widehat
f(n) := \int_{\T} f \overline{e_n} \,d\lambda$.
For $p \ge 1$, let $H^p(\T)$ denote the Hardy space of those $f \in
L^p(\lambda)$ with $\widehat f(n) = 0$ for all $n < 0$.
Write $P_{H^2}$ for the orthogonal projection from $L^2(\lambda) \to
H^2(\T)$.
For all $f \in H^1(\T)$, we have $\GM(|f|) > 0$
(see \ref b.Rudin:RCA/, Theorem 17.17, p.~344).
For the converse, for any function
$f \ge 0$ with $\log f \in L^1(\lambda)$, define
$$
\Phi_f(z) := \exp {1 \over 2} \int_{\T} {e_1(t) + z
\over e_1(t) - z} \log f(t) \,d\lambda(t)
\label e.defPhi_f
$$
for $|z| < 1$.
The {\bf outer} function
$$
\varphi_f(t) := \lim_{r \uparrow 1} \Phi_f(r e_1(t))
\label e.defvarphi_f
$$
exists for $\lambda$-a.e.\ $t \in \T$ and satisfies $|\varphi_f|^2 = f$
$\lambda$-a.e.
Also, $\varphi_f \in H^p(\T)$ iff $f \in L^{p/2}(\lambda)$.
If $f \in L^2(\lambda)$, then the
limit in \ref e.defvarphi_f/ also holds in $L^{1}(\lambda)$.
See \ref b.Rudin:RCA/, Theorem 17.11, p.~340, and Theorem 17.16, p.~343.
Cauchy's integral formula shows that if $f \in L^2(\lambda)$, then
$$
\widehat {\varphi_f}(0) = \Phi_f(0) = \sqrt{\GM(f)}
\,.
\label e.PhiGM
$$
%
By the factorization Theorem 17.17, p.~344, of \ref b.Rudin:RCA/,
if $g \in H^p(\T)$ and $\varphi$ is an outer function such that 
$g \varphi \in L^q(\lambda)$, then $g \varphi \in H^q(\T)$.
Define $\Phi_f := \varphi_f := \constant 0$ if $\log f \notin L^1(\lambda)$.

\bsection {Positive Measures}{s.positive}

Note that the determinant $D_n(\mu)$
is the same when $j, k$ take values in any index set
of $n+1$ consecutive integers.
We shall, in fact, use the index set $ \{ -1, -2, \ldots, -n \}$ for
$D_{n-1}(\mu)$.

In this section, we extend \ref t.first/ to more general determinants
as follows.

\procl t.nonneg
Let $\mu$ be a finite positive measure on $\T$ with infinite support.
Let $f:= [d\mu/d\lambda]$ be the Radon-Nikod\'ym derivative of the
absolutely continuous part of $\mu$.
Given any functions $f_0, \ldots, f_{r}, g_0, \ldots, g_r \in L^2(\mu)$, let
$F_j := P_{H^2}(f_j \overline {\varphi_f})$ and $G_j := P_{H^2}(g_j \overline
{\varphi_f})$.
Define 
$$
p_j := \cases{
f_j &if $0 \le j \le r$,\cr
e_j &if $-1 \ge j \ge -n$\cr
}
$$
and
$$
q_j := \cases{
g_j &if $0 \le j \le r$,\cr
e_j &if $-1 \ge j \ge -n$.\cr
}
$$
We have 
$$
\lim_{n \to\infty}
D_{n-1}(\mu)^{-1}
\det\left[\int p_j \overline{q_k} \,d\mu\right]_{-n \le j, k \le r}
=
\det\left[ \int F_j \overline{G_k} \,d\lambda\right]_{0 \le j, k \le r}
\,.
$$
\endprocl

\procl r.defined
If $\GM(f) = 0$, then $\varphi_f = \constant 0$, so
$F_j = G_j = \constant 0$. Otherwise, $f > 0$ $\lambda$-a.e., so that
$\lambda \ll f\lambda$ and there is no ambiguity about the equivalence
class of $f_j$ or $g_j$ with respect to $\lambda$.
Also, $\int |f_j \overline{\varphi_f}|^2 \,d\lambda = \int |f_j|^2 f \,d\lambda
\le \int |f_j|^2 \,d\mu < \infty$, so that $f_j \overline{\varphi_f} \in
L^2(\lambda)$ and $F_j$ is well defined. Likewise $G_j$ is well defined.
\endprocl

\procl r.nonsing
Since $\mu$ has infinite support, $[\widehat\mu(k-j)]_{-1 \ge j, k \ge -n}$
is non-singular. Indeed, if it were singular, then since it is a Gram
matrix $[\ip{e_j, e_k}_\mu]$, where the subscript $\mu$ indicates that
the inner product is taken in $L^2(\mu)$, it would follow
that the vectors $e_{-1}, \ldots, e_{-n}$ would be linearly dependent,
i.e., there would be scalars $a_j$ such that $\sum_j a_j e_j = 0$
$\mu$-a.e. This would imply that $\mu$ would have support contained in the
zero set of this trigonometric polynomial, i.e., $\mu$ would have support
of cardinality at most $n-1$. 
\endprocl

\procl r.bdtw
The case considered by \ref b.BD/ and \ref b.TW/ is that where all
functions $f_j$ and $g_j$ are of the form $e_n$ for various $n \ge 0$.
These give minors of the Toeplitz matrix other than merely $D_n(\mu)$.
In this case, when $f_j := e_j$ and $g_k := e_k$,
the limiting matrix entries $\int F_j \overline{G_k} \,d\lambda$
become 
$$
\int P_{H^2}(e_j \overline{\varphi_f}) \overline{P_{H^2}(e_k
\overline{\varphi_f})}\,d\lambda 
=
\sum_{l=0}^{\min(j, k)} \overline{\widehat\varphi_f(j-l)}
\widehat\varphi_f(k-l)
\,.
\label e.twform
$$
\refbauthor{BD} gave a different formula than \ref e.twform/; \refbauthor{TW}
gave the same formula as ours.
Both sets of authors assumed that $\mu$ was absolutely continuous.
In addition, \refbauthor{BD} assumed that $f = e^g$ for some $g$ satisfying
$\sum_{n \in \Z} \big( |\widehat g(n)| + |n \widehat g(n)|^2 \big) < \infty$,
while \refbauthor{TW} assumed that $f$ was bounded above and bounded away
from 0.
On the other hand, \refbauthor{BD} showed the strong Szeg\H{o} limit
theorem, which gives finer asymptotics.
\endprocl

\procl r.gren-ros
The special case $r:=0$, $f_0 := g_0 := e_j$ for any fixed $j > 0$ and
$\mu$ is absolutely continuous is due
to Kolmogorov and Wiener (see \ref b.GZ/, Section 10.9).
\endprocl

\procl r.noneed
In case one of $F_j$ or $G_k$ is easier to calculate than the other, one
could use instead of $\int F_j \overline{G_k} \,d\lambda$ either of the
equivalent expressions
$\int F_j \overline{g_k} \varphi_f \,d\lambda$ or
$\int f_j \overline{G_k \varphi_f} \,d\lambda$.
\endprocl

\proofof t.nonneg
For the ease of the reader, we treat first the case $r=0$, $f_0 = g_0 =
\constant 1$, when \ref t.nonneg/ becomes the Szeg\H{o} limit theorem.
Since $\widehat \mu(k-j) = \ip{e_j, e_k}_\mu$, we have that 
$$
D_n(\mu)/D_{n-1}(\mu)
=
\|P_n e_0\|_\mu^2
\,,
$$
where $P_n$ is the orthogonal projection onto $\{e_{-1}, \ldots,
e_{-n}\}^\perp$ in $L^2(\mu)$.
(This is sometimes called ``Gram's formula".)
This quotient therefore tends (monotonically)
to $\|P_\infty e_0\|_\mu^2$, where $P_\infty$ is the
orthogonal projection of $L^2(\mu)$ onto 
$$
H_\infty := \{ g \in L^2(\mu) \st \all {n < 0} \ip{g, e_n}_\mu = 0 \}
\,.
$$
Now $g \in H_\infty$ iff $g \in L^2(\mu)$ and $g \mu$ is an analytic measure.
By the F.\ and M.\ Riesz theorem (\ref b.Rudin:RCA/, Theorem 17.13,
p.~341), it follows that $g \in H_\infty$ iff $g = 0$ a.e.\ with respect to
the singular part of $\mu$, $g \in L^2(f)$ and $g f \in H^1(\T)$.
In particular, we may from now on disregard the singular part of $\mu$.
That is, $P_\infty e_0$ is the same as the orthogonal projection of $e_0$ in
$L^2(f)$ onto $H_\infty := \{ g \in L^2(f) \st \all {n < 0} \ip{g, e_n}_f = 0
\}$ and its norm in $L^2(\mu)$ is the same as its norm in $L^2(f)$.
Write $h_0 := f \cdot P_\infty e_0$
and $\varphi := \varphi_f$.

If $\GM(f) > 0$, then
$h_0/f \in L^2(f)$, $h_0 \in H^1(\T)$ and $h_0/\varphi \in H^2(\T)$.
Also, for all $g \in H_\infty$, we have 
$$
\bigip{h_0/f, g}_\mu = \ip{e_0, g}_\mu
\,.
$$
For all $m \ge 0$, we have $e_m/\overline\varphi \in H_\infty$, whence
$$
\bigip{h_0/f, e_m/\overline\varphi}_\mu = \ip{e_0, e_m/\overline\varphi}_\mu
\,,
$$
or in other words, 
$
\widehat{(h_0/\varphi)}(m)
=
\widehat{\overline\varphi}(m)
$.
Since $\widehat{\overline\varphi}(m) = \sqrt{\GM(f)} \delta_{0, m}$ for $m \ge
0$, we obtain that
$$
\widehat{(h_0/\varphi)}(m)
=
\sqrt{\GM(f)} \delta_{0, m}
\,.
\label e.FC
$$
for $m \ge 0$. 
Since $h_0/\varphi \in H^2(\T)$, \ref e.FC/ holds for all $m$.
That is, $h_0 = \sqrt{\GM(f)} \varphi$.
Therefore, 
$$
\|P_\infty e_0\|_\mu^2 
=
\|h_0/f\|_\mu^2
=
\int {|h_0|^2 \over f} \,d\lambda
=
\GM(f)
\,.
$$

If $\GM(f) = 0$, then for all $g \in H_\infty$, 
$$
\GM(|g f|)^2
=
\GM(|g f|^2)
=
\GM(|g|^2 f) \GM(f) 
=
0
$$
since $\GM(|g|^2 f) \le \int |g|^2 f \,d\lambda < \infty$ as $g \in L^2(f)$.
Since $g f \in H^1(\T)$, this means that $g f = \constant 0$ as noted in \ref
s.notn/. In other words, $g = 0$ $\mu$-a.e.
Therefore $H_\infty = 0$ and so the limit is 0.

We have thus proved the Szeg\H{o} formula.
This proof also shows immediately that the linear span of
$ \{ e_n \st n \ge 0 \}$ is dense in $L^2(\mu)$ iff $\GM(f) = 0$, a theorem
of Kolmogorov and Kre\u\i n.
(More precisely, as written, this proof decides the density of the linear span
of $ \{ e_n \st n \le -1 \}$ by deciding whether its orthocomplement is 0, but
this is equivalent.)

Now we continue with the general case.
Consider $j \ge 0$.
Since 
$$
P_n f_j = f_j - \sum_{-1 \ge i \ge -n} a_i e_i
$$
for some constants $a_i$, row operations can be used to change the $j$th
row from its initial value $[\ip{f_j, q_k}_\mu]_{-n \le k \le r}$ to $[\ip{
P_n f_j, q_k}_\mu]_{-n \le k \le r}$ without changing the determinant.
Since $P_n$ is an orthogonal projection, we have $\ip{P_n f_j, q_k}_\mu =
\ip{P_n f_j, P_n q_k}_\mu$.
If we change all rows $j \ge 0$ in this manner, we obtain a block diagonal
matrix, which shows that
$$
D_{n-1}(\mu)^{-1}
\det\left[\ip{p_j, q_k}_\mu\right]_{-n \le j, k \le r}
=
\det[\ip{P_n f_j, P_n g_k}_\mu]_{0 \le j, k \le r}
\,.
$$
Thus, the limit is
$$
\det[\ip{P_\infty f_j, P_\infty g_k}_\mu]_{0 \le j, k \le r}
\,.
$$
As before, if $\GM(f) = 0$, then $H_\infty = 0$ and the limit is 0. Otherwise, 
the reasoning that led to \ref e.FC/ now leads to 
$$
[f (P_\infty f_j)/\varphi]^{\widehat{\quad}}(m) = \widehat{f_j \overline\varphi}(m)
$$
for all $m \ge 0$, whence $f (P_\infty f_j)/\varphi = F_j$.
Likewise, $f (P_\infty g_k)/\varphi = G_k$.
This gives the formula since 
$$
\ip{P_\infty f_j, P_\infty g_k}_\mu 
=
\int P_\infty f_j \cdot \overline{P_\infty g_k} \,d\mu
=
\int P_\infty f_j \cdot \overline{P_\infty g_k} \cdot f\,d\lambda
=
\int F_j \overline{G_k} \,d\lambda
\,.
\Qed
$$

\procl r.genfn
A bivariate generating function for the matrix entries of \ref e.twform/ is 
$$
\sum_{j, k \ge 0} z^j \overline\zeta^k 
\sum_{l=0}^{\min(j, k)} \overline{\widehat\varphi_f(j-l)}
\widehat\varphi_f(k-l)
=
{\Phi_f(z) \overline{\Phi_f(\zeta)} \over 1 - \overline\zeta z}
\,.
$$
\endprocl

\bsection{Complex Measures}{s.cplx}

We now consider the case of absolutely continuous complex measures, $\mu$.
The proof of the main result in this section,
\ref t.cplx/, could be modified so as to allow a positive
singular part to $\mu$ and to include all of \ref t.nonneg/.
However, the proof would become less elegant.

Let $\trig_n$ denote the linear span of $ \{ e_0, e_1, \ldots, e_n \}$.
Given a pair of functions $\varphi, \psi \in L^2(\lambda)$, consider the
condition 
$$
\texists {\epsilon > 0} \texists {n_0}
\all {n \ge n_0} \all {S \in \trig_n} \texists {T \in \trig_n \setminus \{
0 \}}  \quad
\ip{\varphi S, \psi T}_\lambda
\ge
\epsilon \|\varphi S\|_\lambda \|\psi T\|_\lambda
\,.
\label e.good
$$
Of course, this holds if $\varphi = \psi$, since we may then take $\epsilon
:= 1$ and $T := S$.
Some readers may prefer the following restatement of \ref e.good/. Given
two subspaces $H_1$ and $K_1$ of a Hilbert space $H$, define 
$$
\epsilon(H_1, K_1; H) := \epsilon(H_1, K_1) := \inf_{{x \in H_1 \atop \|x\|
= 1}} \sup_{{y \in K_1 \atop \|y\|=1}} |(x, y)|
\,.
$$
The condition $\epsilon(H_1, K_1) > 0$ is weaker than $H_1 = K_1$ and
stronger than $H_1 \cap K_1^\perp = 0$.
Our condition \ref e.good/ is equivalent to
$$
\liminf_{n \to\infty} \epsilon\big(\varphi \cdot \trig_n, \psi \cdot \trig_n;
L^2(\lambda)\big) > 0
\,.
$$

Condition \ref e.good/ will be used via the following criterion.
We write $H_n \uparrow H_\infty$ to mean that $H_n \subseteq H_{n+1}$ for all
$n$ and $\bigcup H_n$ is dense in $H_\infty$.

\procl l.projection
Suppose that $H$ is a Hilbert space, $H_n, K_n$ are non-zero
closed subspaces for $1 \le n \le \infty$
with $H = H_n + K_n^\perp$ and
$H_n \cap K_n^\perp = 0$ for all $1 \le n \le \infty$.
Suppose that $H_n \uparrow H_\infty$ and $K_n \uparrow K_\infty$.
Let $T_n : H \to K_n^\perp$ be the linear projection along $H_n$ ($1 \le n \le
\infty$).
Then $T_n \to T_\infty$ in the strong operator topology iff
$$
\liminf_{n \to\infty} \epsilon(H_n, K_n) > 0
\,.
\label e.gengood
$$
\endprocl

This lemma should be known, but we could not locate a reference.
Thus, we include its proof in an appendix.
Note that when $H_n = K_n$, which will correspond to the case $\varphi =
\psi$ in our application, it is trivial that $T_n \to T_\infty$ in the
strong operator topology.

As we have noted already, $\epsilon(H_n, K_n) > 0$
implies that $H_n \cap K_n^\perp = 0$.
If $\dim H_n = \dim K_n <\infty$, as will be the case in our application of
\ref e.gengood/, this in turn implies that $H = H_n + K_n^\perp$.

\procl t.cplx
Suppose that $\mu = \psi\overline{\varphi }\lambda$ for some pair of outer
functions $\varphi, \psi \in H^2(\T)$ that satisfies condition \ref e.good/.
Given any functions $f_0, \ldots, f_{r}, g_0, \ldots, g_r \in L^2(|\varphi|^2
+ |\psi|^2)$, let $F_j := P_{H^2}(f_j \overline \varphi)$ and $G_j :=
P_{H^2}(g_j \overline \psi)$.
Define 
$$
p_j := \cases{
f_j &if $0 \le j \le r$,\cr
e_j &if $-1 \ge j \ge -n$\cr
}
$$
and
$$
q_j := \cases{
g_j &if $0 \le j \le r$,\cr
e_j &if $-1 \ge j \ge -n$.\cr
}
$$
We have 
$$
\lim_{n \to\infty}
D_{n-1}(\mu)^{-1}
\det\left[\int p_j \overline{q_k} \,d\mu\right]_{-n \le j, k \le r}
=
\det\left[ \int F_j \overline{G_k} \,d\lambda\right]_{0 \le j, k \le r}
\,.
$$
\endprocl

Note that $L^2(|\varphi|^2 + |\psi|^2) \subseteq L^2(|\mu|)$ by the
Cauchy-Schwarz inequality.

\proof
Let $H_n(\varphi) := e_1 \varphi \trig_n$.
By virtue of \ref e.good/, we have for $n \ge n_0$,
$$
\overline{H_n(\varphi)} \cap \overline{H_n(\psi)}^\perp = 0
\,,
$$
and so
$$
L^2(\T) = \overline{H_n(\varphi)} + \overline{H_n(\psi)}^\perp
\,.
$$
A consequence of Beurling's theorem (\ref b.Rudin:RCA/, Theorem 17.23,
p.~350) is that 
$$
H_n(\varphi) \uparrow H^2_0(\T) := e_1 H^2(\T)
\,.
\label e.beurling
$$
Thus the projection along $\overline{H_n(\varphi)}$ to
$\overline{H_n(\psi)}^\perp$ tends to the orthogonal projection
$P_{\overline{H^2_0(\T)}^\perp} = P_{H^2}$.

Now $\int p_j \overline{q_k} \,d\mu = \ip{\overline\varphi p_j, \overline\psi 
q_k}_\lambda$.
Let $F^{(n)}_j$ be the projection of $\overline\varphi f_k$ along
$\overline{H_n(\varphi)}$ to $\overline{H_n(\psi)}^\perp$.
Row operations show that for $n \ge n_0$, 
$$
D_{n-1}(\mu)^{-1}
\det[\int p_j \overline{q_k} \,d\mu]_{-n \le j, k \le r}
=
\det[\ip{F^{(n)}_j, \overline\psi g_k}_\lambda]_{0 \le j, k \le r}
\,.
$$
Because of our assumption \ref e.good/ and \ref l.projection/, the limit is
$
\det[\ip{F_j, \overline\psi g_k}_\lambda]_{0 \le j, k \le r}
$,
which is the same as
$
\det[\ip{F_j, G_k}_\lambda]_{0 \le j, k \le r}
$.
\Qed

\procl r.beurling
The limit \ref e.beurling/ is often used to prove Beurling's theorem and it
has a simple direct proof: If $g \in H^2_0(\T)$ and $g \perp H_n(\varphi)$
for all $n \ge 0$, then $\widehat{g \overline\varphi}(k) = 0$ for all $k \ge
1$, i.e., $g \overline\varphi \in \overline{H^1(\T)}$.
Dividing by $\overline\varphi$, we get that $g \in \overline{H^2(\T)} =
(H^2_0(\T))^\perp$, so that $g = \constant 0$.
\endprocl

The case considered by \ref b.BD/ and \ref b.TW/ is that where all
functions $f_j$ and $g_j$ are of the form $e_n$ for various $n \ge 0$.
In addition, \refbauthor{BD} assumed that $f = e^g$ for some $g$ satisfying
$\sum_{n \in \Z} \big( |\widehat g(n)| + |n \widehat g(n)|^2 \big) < \infty$, 
while \refbauthor{TW} assumed that
$\varphi$ and $\psi$ are bounded above and that the Toeplitz matrix
corresponding to $\psi\overline\varphi$ has uniformly invertible finite
sections.

The assumption of \refbauthor{BD} implies that of \refbauthor{TW}.
Indeed, write $g = g_1 + g_2$, where $g_1 := \sum_{n \ge 0} \widehat
g(n) e_n$ and $g_2 := \sum_{n < 0} \widehat g(n) e_n$.
Set $f_1 := e^{g_1 + \overline{g_1}}$ and
$f_2 := e^{g_2 + \overline{g_2}}$.
Then $f = \psi \overline\varphi$ with $\psi := e^{g_1} = \varphi_{f_1}$
and $\varphi := e^{\overline{g_2}} = \varphi_{f_2}$, so that $\psi$ and
$\varphi$ are bounded outer functions.
Furthermore, since
$g$ is continuous, Kre\u\i{}n's Theorem (\ref
b.BS:large/, Theorem 1.15, p.~18) in combination with a theorem of Gohberg
and Feldman (\ref b.BS:large/, Theorem 2.11, p.~39) shows that the Toeplitz
matrix of $f$ has uniformly invertible finite sections.

Our theorem covers that of \refbauthor{TW} since boundedness of $\psi$ and
uniform invertibility of finite
sections implies uniform boundedness of the projections along
$\overline{H_n(\varphi)}$ to $\overline{H_n(\psi)}^\perp$, as we see by
simply writing the equations: If $g \in L^2(\lambda)$ is written as $g = u + v$
with $u \in \overline{H_n(\varphi)}$ and $v \in \overline{H_n(\psi)}^\perp$,
then write $\overline u = \sum_{k=1}^n a_k e_k \varphi$. The coefficients
$a_k$ are determined by the requirement that $g - u \perp
\overline{H_n(\psi)}$, i.e., by the equations 
$$
\all {k \in [1, n]} \quad
\widehat{\psi g}(-k)
=
\sum_{j=1}^n {a_j} \widehat{\psi\overline\varphi}(j-k)
\,.
$$
Now observe that $\left[\sum_{k=1}^n |\widehat{\psi g}(-k)|^2 \right]^{1/2}
\le \|\psi g\|_\lambda \le \|\psi\|_\infty \|g\|_\lambda$.

We next give some additional cases when \ref e.good/ holds.
We begin with a reformulation of \ref e.good/, for which we are grateful to
Doron Lubinsky.
Let $w := |\psi|^2$ and $\sigma := \varphi/\psi$.
Then 
$$
\ip{\varphi S, \psi T}_\lambda
=
\int \sigma S \overline T w\,d\lambda
\,,
$$
so that 
$$
\sup_{T \in \trig_n} 
{|\ip{\varphi S, \psi T}_\lambda| \over \|\psi T\|_\lambda}
=
\sup_{T \in \trig_n} 
{|\ip{\sigma S, T}_w| \over \|T\|_w}
=
\|P_{\trig_n}(\sigma S)\|_w
\,,
$$
where the orthogonal projection onto $\trig_n$ takes place in $L^2(w)$.
Note that $\sigma \in L^2(w)$ since $\varphi \in L^2(\lambda)$.
Let $\onp_n$ be the standard orthogonal polynomials for the weight $w$, i.e.,
$\onp_n \in \trig_n$ with positive leading coefficient and $\ip{\onp_m,
\onp_n}_w = \delta_{m, n}$.
If we write 
$$
\sigma S = \sum_{k \ge 0} a_k \onp_k
\,,
\label e.defa_k
$$
then $\|\sigma S\|^2_w =
\sum_{k \ge 0} |a_k|^2$ and
$\|P_{\trig_n}(\sigma S)\|_w^2 = \sum_{k=0}^n |a_k|^2$.
Thus, \ref e.good/ is equivalent to
$$
\texists {\epsilon > 0} \texists {n_0}
\all {n \ge n_0} \all {S \in \trig_n} 
\sum_{k=0}^n |a_k|^2
\ge
\epsilon 
\sum_{k=0}^\infty |a_k|^2
\,,
\label e.goodonp
$$
where $a_k$ are defined by \ref e.defa_k/.

\procl p.poly
If $\varphi, \psi \in L^2(\T)$ and $\sigma := \varphi/\psi$
is an analytic polynomial that has 
no zeroes in the closed unit disc, then \ref e.good/ holds.
\endprocl

Of course, this means that \ref t.cplx/ applies to the pair $\varphi, \psi$
if, in addition, they are outer functions.

\proof 
We use the notation above and show that \ref e.goodonp/ holds.
Let the zeroes of $\sigma$ be $z_1, \ldots, z_m$, all outside the
closed unit disc, and suppose first that each zero is simple.
Let $S \in \trig_n$.
Since $(\sigma S)(z_j) = 0$, we have 
$$
\sum_{l=1}^m a_{n+l} \onp_{n+l}(z_j)
=
-\sum_{k=0}^n a_{k} \onp_k(z_j)
$$
for $1 \le j \le m$ in the notation of \ref e.defa_k/.
Write these equations as 
$$
\sum_{l=1}^m a_{n+l} \tilonp_{n+l}(z_j)/z_j^{n+1}
=
-\sum_{k=0}^n a_{k} \tilonp_k(z_j)/z_j^{n+1}
=: \zeta_j
\,,
\label e.system
$$
where 
$$
\tilonp_k(z) :=
\onp_k(z) \left(\overline{\Phi_w(\overline z^{-1})}\right)^{-1}
\,.
$$
Recall Szeg\H{o}'s asymptotics 
$$
\hbox{for all } \rho > 1,\quad \lim_{k \to\infty} z^{-k} \tilonp_k(z)
=
1 \hbox{ uniformly for } |z| > \rho
\label e.asymp
$$
(see \ref b.GZ/, p.~51).
Because $|z_j| > 1$, it follows that 
$$
C := \sup_{j, n} \sum_{k=0}^n |\tilonp_k(z_j)/z_j^{n+1}|^2 <\infty
\,.
$$ 
By the Cauchy-Schwarz inequality, we obtain that
$$
|\zeta_j|^2 \le C \sum_{k=0}^n |a_k|^2
\,.
$$
Let $M_n := [\tilonp_{n+l}(z_j)/z_j^{n+1}]_{1 \le j, l \le m}$ be the matrix
of coefficients in the system of equations \ref e.system/ for $a_{n+l}$,
considered as variables.
If $M_n$ is not singular, then 
$$
[a_{n+l}]_{1 \le l \le m} = M_n^{-1} [\zeta_j]_{1 \le j \le m}
\,.
$$
Now by \ref e.asymp/, the matrix $M_n$ tends to the Vandermonde
matrix determined by $z_1, \ldots, z_m$.
Hence, for all large $n$, we have not only that $M_n$ is nonsingular, but
also that its inverse has $\ell^2$-norm bounded by some constant $D$
that depends only on $\sigma$.
Therefore, for all large $n$, we have 
$$
\sum_{l=0}^m |a_{n+l}|^2
\le
D^2 \sum_{j=1}^m |\zeta_j|^2
\le
C D^2 m \sum_{k=0}^n |a_k|^2
\,.
$$
This clearly implies \ref e.goodonp/ for $\epsilon := 1/(1+C D^2 m)$
and finishes the proof for the case of simple zeroes.

Now suppose that the zero $z_j$ of $\sigma$ has multiplicity $r_j$, so that
$s := \sum_{j=1}^m r_j$ is the degree of $\sigma$.
Multiply \ref e.asymp/ by $z^p$ and
take the $r$th derivative (\ref b.Rudin:RCA/, Theorem 10.28, p.~214)
to obtain that
$$
\lim_{k \to\infty} 
\qonp_{k, p}^{(r)}(z) 
=
z^{p-r} \prod_{t=0}^{r-1} (p-t)
\quad\hbox{ for } |z| > 1
\,,
\label e.asymp-der
$$
where
$$
\qonp_{k, p}^{(r)}(z) :=
\left({d \over dz}\right)^r \left( z^{-k+p} \tilonp_k(z) \right)
\,.
$$
Since $z \mapsto z^{-n-1} \sigma(z) S(z)
\left(\overline{\Phi_w(\overline z^{-1})}\right)^{-1}$
has a zero at $z_j$ of order at least $r_j$, it follows that
$$
\sum_{l=1}^s a_{n+l} \qonp^{(r)}_{n+l, l-1}(z_j)
=
-\sum_{k=0}^n a_{k} \qonp^{(r)}_{k, l-1}(z_j)
$$
for $0 \le r < r_j$ and $1 \le j \le m$.
We may now follow the same reasoning as for the case of simple zeroes, but
instead of finding a coefficient matrix tending to a Vandermonde matrix,
we find instead a limit matrix that has $r_j$ columns corresponding to each
$z_j$, namely, for each $r = 0, 1, \ldots, r_j - 1$, it has the column
$$
\left[z_j^{l-1-r} \prod_{t=0}^{r-1} (l-1-t) \right]_{1 \le l \le s}
\,.
$$
Thus, by reasoning analogous to before, it suffices to establish that these
columns form a nonsingular matrix.
To do this, we show that there is no nontrivial linear relation among the
rows.
Indeed, if $b_1, \ldots, b_s$ are constants such that for $1 \le j \le m$
and $0 \le r < r_j$, 
$$
\sum_{l=1}^s b_l z_j^{l-1-r} \prod_{t=0}^{r-1} (l-1-t) = 0
\,,
$$
then the polynomial $\sum_{l=1}^s b_l z^{l-1}$ has a zero at $z_j$ of order
at least $r_j$ for each $1 \le j \le m$. But since this polynomial has
degree at most $s-1$, this implies that the polynomial is identically zero,
i.e., all $b_l = 0$, as desired.
\Qed

It would be interesting to have a good characterization of those $\varphi,
\psi$ such that \ref e.good/ holds.

\bsection{Higher Dimensions}{s.high}

The technology we use to replace the theory of Hardy spaces for higher
dimensions was provided by \ref b.HelLow/, who proved the extension of \ref
t.first/ to higher dimensions. We review the relevant definitions and facts
from their theory before giving our theorems, which extend Theorems \briefref
t.nonneg/ and \briefref t.cplx/. Fix a positive integer
$d$ and let $\lambda$ be Lebesgue measure on $\T^d := \R^d/\Z^d$.

For $k \in \Z^d$ and $x \in \T^d$, let $e_k(x) := e^{2\pi i k \cdot x}$.
For $f \in L^2(\T^d)$, write $\widehat f(k) := (f, e_k)_\lambda := \int_{\T^d}
f \overline{e_k} \,d\lambda$.
Let $S \subset \Z^d$ have the properties $S \cup (-S) = \Z^d \setminus \{ \bfz
\}$, $S \cap (-S) = \emptyset$, and $S+S \subset S$. The associated {\bf
ordering} of $\Z^d$ is that where $k \prec l$ iff $l - k \in S$.
For example, we could have $(k_1, k_2, \ldots, k_d) \prec \bfz$
if $k_i < 0$ when $i$ is the first index such that $k_i \ne 0$, which we
call the {\bf lexicographic ordering}.
The replacement for the Hardy spaces $H^p(\T)$ ($p \ge 1$) are the {\bf
Helson-Lowdenslager spaces}
$$
\hl^p := \hl^p(\T^d, S) := \Big \{ \varphi \in L^p(\T^d) \st \supp \widehat
\varphi \subset S \cup \{ \bfz \} \Big \}
\,.
$$
Let $P_{\hl^2} : L^2(\T^d) \to \hl^2$ be the orthogonal projection $\sum_{k
\in \Z^d} a_k e_k \mapsto \sum_{k \in S \cup \{ \bfz \}} a_k e_k$.
For $0 \le f \in L^1(\lambda)$ and any set $R \subseteq \Z^d$,
let $[R]$ be the linear span of $ \{ e_k \st k \in R\}$ and $[R]_f$ be
its closure in $L^2(f)$.
In place of outer functions, we use {\bf spectral factors},
which are the functions $\varphi \in \hl^2$ with the properties 
$$
\widehat\varphi(\bfz) > 0
\label e.at0
$$
and
$$
1/\varphi \in [S \cup \{ \bfz \}]_{|\varphi|^2}
\,.
\label e.inverse
$$
\ref b.HelLow/ show that for $0 \le f \in L^1(\T^d)$, the condition $\GM(f) >
0$ is equivalent to the existence of a spectral factor $\varphi$ such that
$|\varphi|^2 = f$. (More precisely, they prove $\GM(f) > 0$ iff $\exists
\varphi \in \hl^2$ satisfying 
\ref e.at0/. Their proof
shows that in this case, $\varphi$ can be chosen so that also \ref e.inverse/
holds.)

\procl l.divide
Suppose that $\varphi$ satisfies \ref e.inverse/, $h \in \hl^1$, and
$h/\varphi \in L^2(\T^d)$. Then $h/\varphi \in \hl^2$.
\endprocl

\proof
Let $f := |\varphi|^2$. By \ref e.inverse/, there exist trigonometric
polynomials $p_n$ with $\supp \widehat {p_n} \subset S \cup \{ \bfz \}$ and
such that $p_n \to 1/\varphi$ in $L^2(f)$.
Since $h/\varphi \in L^2(\T^d)$, we have $h/f \in L^2(f)$.
Thus, 
\begineqalno
\widehat{h/\varphi}(k)
&=
\int {h \over \varphi} \overline{e_k} \,d\lambda
=
\int {h \overline{e_k}\over f} {1\over \varphi} f \,d\lambda
=
\ip{h\overline{e_k}/f, 1/\overline\varphi}_f
=
\lim_{n \to\infty} \ip{h\overline{e_k}/f,  \overline{p_n}}_f
\cr&=
\lim_{n \to\infty} \int h p_n \overline{e_k}\,d\lambda
=
\lim_{n \to\infty} \widehat{h p_n}(k)
\,.
\cr
\endeqalno 
Since $h p_n \in \hl^1$, we have $\widehat{h p_n}(k) = 0$ for $k \notin S$, 
whence $\widehat{h/\varphi}(k) = 0$ for $k \notin S$.
That is, $h/\varphi \in \hl^2$.
\Qed

For $A \subseteq \Z^d$, let $(A)$ denote the set of corresponding complex
exponentials $ \{ e_k \st k \in A \}$.
For any two finite ordered sets of functions $\F, \G \subset L^2(\mu)$
of the same cardinality, let 
$$
\bigip{\F, \G}_\mu :=
\det\big[\ip{p, q}_\mu\big]_{p \in \F , q \in \G }
\,.
$$
Here, the ordering of the sets $\F, \G$ is used to order the rows and columns
of the matrix whose determinant appears in this equation.
Also, we write $\ip{p, q}_\mu := \int p \overline q \,d\mu$ even if $\mu$
is a complex measure.
Write $\F \Cup \G$ for the set $\F \cup \G$ ordered by concatenating
$\G$ after $\F$.

We are now ready to state and prove our extension of \ref t.nonneg/.

\procl t.nonneg-d
Let $w : \T^d \to \CO{0,\infty}$ be measurable with
$\GM(w) > 0$. Let $\varphi$ be a spectral factor for $w$.
Given any two finite ordered sets of functions $\F, \G \subset L^2(w)$
of the same cardinality, let $\F' := \Seq{ P_{\hl^2}(f \overline \varphi) \st
f \in \F }$ and define $\G'$ likewise.
Let $S_n \subset -S$ be finite ordered sets increasing to $-S$.
We have 
$$
\lim_{n \to\infty}
{\bigip{\F \Cup (S_n) , \G \Cup (S_n)}_w
\over
\bigip{(S_n), (S_n)}_w}
=
\bigip{\F', \G'}_\lambda
\,.
\label e.past
$$
\endprocl

The case $\F = \G = \Seq{\constant 1}$ is the theorem of \ref b.HelLow/.
Actually, there are special considerations in that case that allow
\refbauthor{HelLow} to make the same conclusion when $w = [d\mu/d\lambda]$
and the matrix entries are given by inner products in $L^2(\mu)$.

\proof
Let $P_n$ be the orthogonal projection onto $(S_n)^\perp$ in $L^2(w)$.
Also, let $P_\infty$ be the
orthogonal projection of $L^2(w)$ onto 
$$
H_\infty := (-S)^\perp
=
\{ g \in L^2(w) \st \all {n \in -S} \ip{g, e_n}_w = 0 \}
\,.
$$
Define $\F'_n := \Seq{P_n f \st f \in \F}$ and likewise for $\G'_n$,
$\F'_\infty$, and $\G'_\infty$.
By row and column operations, we have 
$$
{\bigip{\F \Cup (S_n) , \G \Cup (S_n)}_w
\over
\bigip{(S_n), (S_n)}_w}
=
\bigip{\F_n', \G_n'}_w
\,.
$$
This therefore tends to
$\bigip{\F_\infty', \G_\infty'}_w$ (indeed, we have entry-wise convergence
of the corresponding matrices).
Now $g \in H_\infty$ iff $g \in L^2(w)$ and $\supp \widehat{g w} \subseteq
S \cup \{ \bfz \}$.
Thus, 
$$
H_\infty 
=
\{ g \in L^2(w) \st g w \in \hl^1 \}
\,.
$$

Let $f \in L^2(w)$ and
write $h := w \cdot P_\infty f$.
Since $h/w \in H_\infty$, we have $h/w \in L^2(w)$ and $h \in \hl^1$.
{}From the first of these relations, we see that $h/\varphi \in L^2(\T^d)$.
Also, 
$$
\all {g \in H_\infty}\quad  \bigip{h/w, g}_w = \ip{f, g}_w
\,.
\label e.whatyouget
$$
For all $m \in S \cup \{ \bfz \}$, we have $e_m \varphi \in \hl^2 \subset
\hl^1$, so that $e_m /\overline \varphi \in H_\infty$.
Therefore, \ref e.whatyouget/ implies that
$$
\bigip{h/w, e_m/\overline\varphi}_w = \ip{f, e_m/\overline\varphi}_w
\,,
$$
or in other words, 
$
\widehat{(h/\varphi)}(m)
=
\widehat{f\overline\varphi}(m)
$ for all $m \in S \cup \{ \bfz \}$.
Since $h/\varphi \in \hl^2$ by \ref l.divide/, it follows that
$h/\varphi = P_{\hl^2}(f \overline\varphi)$.
Thus, we have proved that for all $f$, we have 
$$
P_\infty f = {P_{\hl^2}(f \overline\varphi) \over \overline\varphi}
\,.
\label e.projform
$$

This gives \ref e.past/ since 
$$
\ip{P_\infty f, P_\infty g}_w 
=
\int P_\infty f \cdot \overline{P_\infty g} \cdot w\,d\lambda
=
\int P_{\hl^2}(f \overline\varphi) \overline{P_{\hl^2}(g \overline\varphi)}
\,d\lambda
\,.
\Qed
$$

Unlike in \ref t.nonneg/, the limit \ref e.past/ is not necessarily 0 when
$\GM(w) = 0$. For example, let $S$ be the lexicographic ordering
on $\Z^2$ and $w(x_1, x_2)$ be a function that depends only on $x_2$ with
$\GM(w) = 0$.  If $\F := \G := \{ e_{(1, 0)} \}$, then the left-hand side of
\ref e.past/ is equal to $\widehat w(\bfz)$, which need not equal 0.

However, if the order is archimedean, such as if $S = \{ k \in \Z^d \st k
\cdot x > 0 \}$, where $x \in \R^d$ has at least two coordinates whose
quotient is irrational, then the limit \ref e.past/ is 0 when $\GM(w) = 0$,
as we now show.
(All archimedean orders arise in this way; in fact, for a characterization of
all orders, see \ref b.Teh/, \ref b.Zaiceva/, or \ref b.Trevisan/.)

\procl p.GM0
Let $w : \T^d \to \CO{0,\infty}$ be measurable with $\GM(w) = 0$. 
Suppose that the order induced by $S$ is archimedean.
Given any two finite ordered sets of functions $\F, \G \subset L^2(\mu)$
of the same cardinality,
and any finite ordered sets $S_n \subset -S$ increasing to $-S$,
we have 
$$
\lim_{n \to\infty}
{\bigip{\F \Cup (S_n) , \G \Cup (S_n)}_w
\over
\bigip{(S_n), (S_n)}_w}
=
0
\,.
$$
\endprocl

\proof
It suffices to establish that $H_\infty = 0$ in the notation of the proof of
\ref t.nonneg-d/.
Now for all $g \in H_\infty$, 
$$
\GM(|g w|)^2
=
\GM(|g w|^2)
=
\GM(|g|^2 w) \GM(w) 
=
0
$$
since $\GM(|g|^2 w) \le \int |g|^2 w \,d\lambda < \infty$ as $g \in L^2(w)$.
Since $g w \in \hl^1$, this means that $g w = \constant 0$ by the main result
of \ref b.arens/. 
Thus, $H_\infty = 0$.
\Qed

\procl r.uniq
Suppose that $\GM(w) > 0$ and that $\varphi$ is a spectral factor for $w$.
By \ref e.projform/, we have $w P_\infty \bfo = \varphi
P_{\hl^2}(\overline\varphi) = \overline{\widehat\varphi(\bfz)} \varphi$, so
that $\varphi$ is uniquely determined by $w$.
It is essentially by this formula that \ref b.HelLow/ proved the existence
of a spectral factor.
This method goes back to \ref b.Szego:outer/.
\endprocl

The extension of \ref t.cplx/ is relatively straightforward:

\procl t.cplx-d
Let $S_n \subset -S$ be finite ordered sets increasing to $-S$.
Let $\mu = \psi\overline\varphi\lambda$ for some pair of spectral factors
$\varphi, \psi$ that satisfy the condition
$$
\liminf_{n \to\infty} \epsilon\big(\varphi \cdot [-S_n], \psi
\cdot [-S_n]; L^2(\lambda)\big) > 0
\,.
\label e.good-d
$$
Given any two finite ordered sets of functions
$\F, \G \subset L^2(|\varphi|^2 + |\psi|^2)$
of the same cardinality, let $\F' := \Seq{ P_{\hl^2}(f \overline \varphi) \st
f \in \F }$ and $\G' := \Seq{ P_{\hl^2}(g \overline \psi) \st
g \in \G }$.
We have 
$$
\lim_{n \to\infty}
{\bigip{\F \Cup (S_n) , \G \Cup (S_n)}_\mu
\over
\bigip{(S_n), (S_n)}_\mu}
=
\bigip{\F', \G'}_\lambda
\,.
\label e.past
$$
\endprocl

\proof
Let $H_n(\varphi) := \{ \varphi e_k \st k \in -S_n \}$.
By virtue of \ref e.good-d/, we have for $n \ge n_0$,
$$
\overline{H_n(\varphi)} \cap \overline{H_n(\psi)}^\perp = 0
\,,
$$
and so
$$
L^2(\T^d) = \overline{H_n(\varphi)} + \overline{H_n(\psi)}^\perp
\,.
$$
We claim that
$$
H_n(\varphi) \uparrow \hl^2_0 := \left(\overline{\hl^2}\right)^\perp
\,.
\label e.beurling-d
$$
Indeed, if $g \in \hl^2_0$ and $g \perp H_n(\varphi)$
for all $n \ge 0$, then $\widehat{g \overline\varphi}(k) = 0$ for all $k \in
S$, i.e., $g \overline\varphi \in \overline{\hl^1}$.
By \ref l.divide/, we may divide by $\overline\varphi$ to obtain that $g \in
\overline{\hl^2} = (\hl^2_0)^\perp$, so that $g = \constant 0$.
This proves \ref e.beurling-d/.

Thus the projection along $\overline{H_n(\varphi)}$ to
$\overline{H_n(\psi)}^\perp$ tends to the orthogonal projection
$P_{\overline{\hl^2_0}^\perp} = P_{\hl^2}$.

Now $\ip{f, g}_\mu = \ip{\overline\varphi f, \overline\psi g}_\lambda$ for any
$f, g \in L^2(|\varphi|^2+|\psi|^2)$.
Let $\F_n$ be the image of $\F$ under
the projection along
$\overline{H_n(\varphi)}$ to $\overline{H_n(\psi)}^\perp$.
Let $\G'' := \{ \overline\psi g \st g \in \G \}$.
Row operations show that for $n \ge n_0$, 
$$
{\bigip{\F \Cup (S_n) , \G \Cup (S_n)}_\mu
\over
\bigip{(S_n), (S_n)}_\mu}
=
\bigip{\F_n, \G''}_\lambda
\,.
$$
Because of our assumption \ref e.good-d/ and \ref l.projection/, the limit is
$
\bigip{\F', \G''}_\lambda
$,
which is the same as
$
\bigip{\F', \G'}_\lambda
$.
\Qed

\bsection{Appendix}{s.app}

In order to prove \ref l.projection/, we first demonstrate the following
lemma.

\procl l.proj1
Suppose that $H$ is a Hilbert space, $H_1$ and $K_1$ are non-zero
closed subspaces, $H = H_1 + K_1^\perp$, and $H_1 \cap K_1^\perp = 0$.
Let $T : H \to K_1^\perp$ be the linear projection along $H_1$
and 
$$
\epsilon := \epsilon(H_1, K_1)
\,.
$$
Then 
$$
\|T\| \le {1 \over 2} + \sqrt{{1 \over 2(1 - \sqrt{1 - \epsilon^2})} + {1
\over 4}}
\label e.bddT
$$
and
$$
\epsilon \ge {1 \over \sqrt{1 + \|T\|^2}}
\,.
\label e.bddepsi
$$
\endprocl

\proof
Let $v \in H$ with $\|v\| = 1$.
Write $w := T v \in K_1^\perp$ and $u := v - w \in H_1$.
Choose $y \in K_1$ such that $\ip{u, y} \ge \epsilon \|u\| \|y\|$ and
$\|y\| = \epsilon \|u\|$.
We have 
\begineqalno
1 
&=
\|v\|^2
=
\|u + w\|^2
=
\|u\|^2 + \|w\|^2 + 2 \Re \ip{u, w}
\cr&=
\|u\|^2 + \|w\|^2 + 2 \Re \ip{u - y, w}
\ge
\|u\|^2 + \|w\|^2 - 2 \|u - y\| \| w\|
\cr&=
\|u\|^2 + \|w\|^2 - 2 \| w\| \left[ \|u\|^2 + \|y\|^2 - 2 \Re\ip{u, y}
\right]^{1/2}
\cr&\ge
\|u\|^2 + \|w\|^2 - 2 \| w\| \left[ \|u\|^2 + \|y\|^2 - 2 \epsilon \|u \|
\|y \| \right]^{1/2}
\cr&=
\|u\|^2 + \|w\|^2 - 2 \| w\| \|u\| \sqrt{1- \epsilon^2}
\cr&=
\left(\|u\| - \|w\|\right)^2 + 2 \| w\| \|u\| (1 - \sqrt{1- \epsilon^2})
\cr&\ge
2 \| w\| \|u\| (1 - \sqrt{1- \epsilon^2})
=
2 \| w\| \|v - w\| (1 - \sqrt{1- \epsilon^2})
\cr&\ge
2 \| w\| (\|w\| - 1) (1 - \sqrt{1- \epsilon^2})
\,.
\endeqalno
Simple algebra shows that this inequality implies 
$$
\|T v\| = \|w\|
\le {1 \over 2} + \sqrt{{1 \over 2(1 - \sqrt{1 - \epsilon^2})} + {1
\over 4}}
\,.
$$
Since $\|v\| = 1$, this is equivalent to \ref e.bddT/.

To prove \ref e.bddepsi/, note that by definition of $\epsilon$, we have
that
$$
\texists {x \in H_1 \setminus \{ 0 \}} \all {y \in
K_1} \quad |\ip{x, y}| \le \epsilon \|x\| \|y\|
\,.
\label e.genbad
$$
Choose such an $x$ with $\|x\|= 1$ and set $y := P_{K_1} x$, the orthogonal
projection of $x$ onto $K_1$.
Then $y \ne 0$ because $H_1 \cap K_1^\perp = 0$.
Write $z := x - y \in K_1^\perp$.
Since $y = x - z$ with $x \in H_1$ and $z \in K_1^\perp$, it follows that
$T y = -z$.
By \ref e.genbad/ applied to $y$, we have 
\begineqalno
\|z\|^2
&=
\|x - y\|^2
=
1 + \|y\|^2 - 2 \Re \ip{x, y}
\cr&\ge
1 + \|y\|^2 - 2 \epsilon \|y\|
\label e.interm
\cr&=
1 + (\|y\| - \epsilon)^2 - \epsilon^2
\cr&\ge
1 - \epsilon^2
\,.
\endeqalno
Furthermore, $\|z\|^2 = 1 - \|y\|^2$, whence comparison to \ref e.interm/
shows that $\|y\| \le \epsilon$.
Therefore, 
$$
\|T\| 
\ge
\|T y\|/\|y\|
=
\|z\|/\|y\|
\ge
\sqrt{1 - \epsilon^2}/\epsilon
\,,
$$
which is equivalent to \ref e.bddepsi/.
\Qed

\proofof l.projection
For each $n$ and any $v \in H_n$, we have $T_m v = 0 = T_\infty v$ for all
$m \ge n$.
Also, for each $v \in K_\infty^\perp$, we have $T_m v = v = T_\infty v$ for
all $m$.
Therefore, $T_n v \to T_\infty v$ for all $v$ belonging to the dense
set $\bigcup_n H_n + K_\infty^\perp$.
It follows by continuity and the principle of uniform boundedness that $T_n
\to T_\infty$ in the strong operator topology iff 
$$
\sup_n \|T_n\| <\infty
\,.
\label e.unifbdd
$$
If \ref e.gengood/ holds, then \ref e.unifbdd/ is a consequence of \ref
e.bddT/, while if \ref e.unifbdd/ holds, then \ref e.gengood/ is a
consequence of \ref e.bddepsi/.
\Qed

\medbreak
\noindent {\bf Acknowledgements.}\enspace
I am grateful to Doron Lubinsky for allowing me to include \ref e.goodonp/
and to Jeff Geronimo for a reference.

\bibfile{\jobname}
\def\noop#1{\relax}
\input \jobname.bbl

\filbreak
\begingroup
\eightpoint\sc
\parindent=0pt\baselineskip=10pt

Department of Mathematics,
Indiana University,
Bloomington, IN 47405-5701
\emailwww{rdlyons@indiana.edu}
{http://php.indiana.edu/\string~rdlyons/}

and

School of Mathematics,
Georgia Institute of Technology,
Atlanta, GA 30332-0160
\email{rdlyons@math.gatech.edu}

\endgroup

\bye